\documentclass{agtart_a}
\pdfoutput=1

\title{Twisted Alexander polynomials detect the unknot}

\author{Daniel S Silver}
\givenname{Daniel S}
\surname{Silver}
\address{Department of Mathematics and Statistics\\
University of South Alabama\\\newline
Mobile, AL 36688-0002\\USA}
\email{silver@jaguar1.usouthal.edu}
\urladdr{}

\author{Susan G Williams}
\givenname{Susan G}
\surname{Williams}
\email{swilliam@jaguar1.usouthal.edu}
\urladdr{}

\volumenumber{6}
\issuenumber{}
\publicationyear{2006}
\papernumber{66}
\startpage{1893}
\endpage{1901}

\doi{}
\MR{}
\Zbl{}


\received{5 June 2006}
\revised{20 August 2006}
\accepted{11 September 2006}
\published{14 November 2006}
\publishedonline{14 November 2006}
\proposed{}
\seconded{}
\corresponding{}
\editor{CPR}
\version{}

\arxivreference{math.GT/0604084}

\AtBeginDocument{\let\bar\wbar\let\tilde\wtilde\let\hat\what}

\makeatletter
\def\cnewtheorem#1[#2]#3{\newtheorem{#1}{#3}[section]
\expandafter\let\csname c@#1\endcsname\c@theorem}
\makeatother

\numberwithin{equation}{section}

\cnewtheorem{prop}[theorem]{Proposition}
\cnewtheorem{lemma}[theorem]{Lemma}
\cnewtheorem{cor}[theorem]{Corollary}
\cnewtheorem{conj}[theorem]{Conjecture}
\newtheorem*{theorem*}{Theorem}

\theoremstyle{definition}
\cnewtheorem{definition}[theorem]{Definition}
\cnewtheorem{example}[theorem]{Example}
\cnewtheorem{remark}[theorem]{Remark}

\def\<{{\langle}}
\def\>{{\rangle}}

\def\L{{\Lambda}}
\def\a{{\alpha}}
\def\b{{\beta}}

\begin{document}

\begin{abstract}
 The group of a nontrivial knot admits a finite permutation
 representation such that the corresponding twisted Alexander
 polynomial is not a unit.
\end{abstract}

\maketitle 
%%%%%%%%%%%%%%%%%%%%%%%%%%%%%% 1. INTRODUCTION %%%%%%%%%%%%%
\section{Introduction}
\label{intro} 

Twisted Alexander polynomials of knots in ${\mathbb S}^3$ were introduced  by X-S Lin in \cite{lin}. They were defined more generally for any finitely presentable group with infinite abelianization by M Wada \cite{wada}. Many papers subsequently appeared on the topic. Notable among them is \cite{kirkliv}, by P Kirk and C Livingston,
placing twisted Alexander polynomials of knots in the classical context of abelian invariants. A slightly more general approach by J Cha \cite{cha} permits coefficients in a Noetherian unique factorization domain.

In Hillman--Livingston--Naik \cite{hln}  two examples are given of Alexander polynomial $1$ hyperbolic knots for which twisted Alexander polynomials provide periodicity obstructions. In each case, a finite representation of the knot group is used to obtain a nontrivial twisted polynomial. Such examples motivate the question: Does the group of any nontrivial knot admit a finite representation such that the resulting twisted Alexander polynomial is not a unit (that is, not equal to $\pm t^i)$?

\begin{theorem*} Let $k\subset {\mathbb S}^3$ be a nontrivial knot. There exists a finite permutation representation such that the corresponding twisted Alexander polynomial $\Delta_\rho (t)$ is not a unit. 
\end{theorem*}

A key ingredient of the proof of the theorem is a recent theorem of M Lackenby \cite{lackenby} which implies that some cyclic cover of ${\Bbb S}^3$ branched over $k$ has a fundamental group with arbitrarily large finite quotients. The quotient map pulls back to a representation of the knot group. A result of J Milnor \cite{milnor} allows us to conclude that for sufficiently large quotients, the associated twisted Alexander polynomial is nontrivial.

We are grateful to Abhijit Champanerkar and Walter Neumann for helpful discussions. We thank the referee for many thoughtful suggestions that improved the exposition of the paper.  Both authors were partially supported by NSF grant DMS-0304971.

\section{Preliminary material}
\vspace{-5pt}
\subsection{Review of twisted Alexander polynomials}
\vspace{-5pt}
\label{review} Let $X$ be a finite CW complex. Its fundamental group $\pi =\pi_1 X$ acts on the left of the universal cover $\tilde X$ by covering transformations. 

Assume that $\epsilon$ is an epimorphism from $\pi$ to an infinite cyclic group $\< t \mid \>$. 
Given a Noetherian unique factorization domain $R$, we identify the group ring $R[\< t \mid \>]$ with the ring of Laurent polynomials $\L=R[t, t^{-1}]$. (Here we will be concerned only with the case $R= {\Bbb Z}$.) 

Assume further that $\pi$ acts on the right of a free $R$--module $V$ of finite rank via a representation $\rho\co  \pi \to GL(V)$. Define a $\L$--$R[\pi]$ bimodule structure on $\L \otimes_R V$ by  $t^j (t^n \otimes v) = t^{n+j} \otimes v$ and $(t^n \otimes v) g = t^{n+ \epsilon(g)} \otimes v \rho(g)$ for $v \in V$ and $g \in \pi$. The groups  of the cellular chain complex $C_*(\tilde X; R)$ are left $R[\pi]$--modules. The twisted complex of $X$ is defined to be the chain complex of left $\L$--modules:
$$C_*(X; \L \otimes V) = (\L \otimes V)\otimes_{R[\pi]} C_*(\tilde X; R).$$
The twisted homology $H_*(X; \L\otimes V)$ is the homology of 
$C_*(X; \L \otimes V)$. 

Since $V$ is finitely generated and $R$ is Noetherian, $H_*(X; \L\otimes V)$ is a finitely presentable $\L$--module. Elementary ideals and characteristic polynomials are defined in the usual way. Begin with an $n \times m$ presentation matrix corresponding to a presentation for $H_1(X; \L\otimes V)$ with $n$ generators and $m\ge n$ relators. The ideal in $\L$ generated by the $n \times n$ minors is an invariant of $H_1(X; \L\otimes V)$. The greatest common divisor of the minors, the \emph{twisted Alexander polynomial} of $X$, is an invariant as well. It is well defined up to a unit of $\L$. Additional details can be found in \cite{cha}. An alternative, group-theoretical approach can be found in \cite{sw06}.

In what follows, $X$ will denote the exterior of a nontrivial knot $k$, that is, the closure of ${\Bbb S}^3$ minus a regular neighborhood of $k$. 

\subsection{Periodic representations}  The knot group $\pi$ is a semidirect product $\< x \mid \> \ltimes \pi'$, where $x$ is a meridional generator and $\pi'$ denotes the commutator subgroup $[\pi, \pi]$. Every element has a unique expression of the form $x^j w$, where $j \in {\mathbb Z}$ and $w \in \pi'$.  

For any positive integer $r$, the fundamental group of the $r$--fold cyclic cover $X_r$ of $X$ is isomorphic to $\<x^r \mid \> \ltimes \pi'$. The fundamental group of the $r$--fold cyclic cover $M_r$ of ${\Bbb S}^3$ branched over $k$ is the quotient group $\pi_1(X_r)/ \<\< x^r\>\>$, where $\<\< \cdot \>\>$ denotes the normal closure. Consequently, $\pi_1M_r = \pi'/[\pi', x^r]$. 

\begin{definition}\label{defper} A representation $p\co  \pi' \to \Sigma$ is {\it periodic} with {\it period} $r$ if it factors through $\pi_1 M_r$. If $r_0$ is the smallest such positive number, then $p$ has \emph{least} least period $r_0$. \label{per} \end{definition}

\begin{remark} The condition that $p$ factors through $\pi_1 M_r$ is equivalent  to the condition that $p(x^{-r}a x^r) =p(a)$ for every $a \in \pi'$. \end{remark}

The following is a consequence of the fact that $M_1$ is ${\mathbb S}^3$. 

\begin{prop}\label{period1} If $p\co  \pi'\to \Sigma$ has period $1$, then $p$ is trivial. 
\end{prop}

Assume that $p\co  \pi'\to \Sigma$ is surjective and has least period $r_0$.
We extend $p$ to a homomorphism $P \co  \pi \to \<\xi \mid \xi^{r_0}\>\ltimes_\theta \Sigma^{r_0}$,  mapping $x \mapsto \xi$ and elements $u \in \pi'$ to $(p(u), p(x^{-1}ux), \ldots, p(x^{-({r_0}-1)}ux^{{r_0}-1})) \in \Sigma^{r_0}$. Conjugation by $\xi$ in the semidirect product induces $\theta\co  \Sigma^{r_0} \to \Sigma^{r_0}$ described by $(\alpha_1, \ldots, \alpha_{r_0}) \mapsto (\alpha_2, \ldots, \alpha_{r_0}, \alpha_1)$. The lemma below assures us that the image of $\pi'$ under $P$ has order no less than the order of $p(\pi')$. 

\begin{lemma}\label{order} $|P(\pi')|\ge |p(\pi')|$ \end{lemma}

\begin{proof} The image $P(\pi')$ is contained in $\Sigma^{r_0}$. First coordinate projection $\Sigma^{r_0} \to \Sigma$ obviously maps $P(\pi')$ onto $p(\pi')$. \end{proof}

In what follows we will assume that $\Sigma$ is finite. Hence $P (\pi)$ is also finite, and it is isomorphic to a group of permutations of a finite set acting transitively (that is, the the orbit of any element under $P (\pi)$ is the entire set.) We can ensure that the subgroup $P(\pi')$ also acts transitively, as the next lemma shows. 

We denote the symmetric group on a set ${\cal A}$ by $S_{\cal A}$.

\begin{lemma}\label{permrep} The group $P (\pi)$ embeds in the symmetric group $S_{P (\pi')}$  in such a way that $P(\pi')$ acts transitively. \end{lemma}

\begin{proof} Embed $P(\pi')$ in $S_{P (\pi')}$ via the right regular representation $\psi\co  P (\pi') \to S_{P (\pi')}$.  Given
$ \b = (\b_1, \ldots, \b_{r_0})\in P(\pi')$, the permutation $\psi(\b)$ maps
$(\a_1, \ldots, \a_{r_0})\in P(\pi')$ to $(\a_1 \b_1, \ldots, \a_{r_0} \b_{r_0})$. 
Extend $\psi$ to $\Psi\co  P(\pi) \to S_{P (\pi')}$ by assigning to $\xi$ the permutation of $P(\pi')$ given by $(\a_1, \ldots, \a_{r_0}) \Psi(\xi) = (\a_2, \ldots, \a_{r_0}, \a_1)$. It is straightforward to check that $\Psi$ respects the action $\theta$ of the semidirect product, and hence is a well-defined homomorphism. 

To see that $\Psi$ is faithful, suppose that $\Psi(\xi^i  \b)$ is trivial for some $1\le i < r_0$, $ \b \in P(\pi')$. Then $\Psi(\xi^i) = \psi( \b^{-1})$. By considering the effect of the permutation  on $  1= (1, \ldots, 1)$, we find that $ \b$ must be $ 1$ and hence the action of $\Psi(\xi^i)$ is trivial. It follows that $p$ has period $i <r_0$, contradicting the assumption that $r_0$ is the least period. 
 \end{proof}
 
We summarize the above construction. 
\begin{lemma}\label{summary} Given a finite representation $p \co  \pi' \to \Sigma$ of period $r$, there is a finite permutation representation $P \co  \pi_1 X \to S_N$  such that $P  \vert_{\pi'}$ is $r$--periodic and transitive. Moreover, $|P(\pi')|=N \ge |p(\pi')|$.  \end{lemma}

\vspace{-5pt}
\subsection{Twisted Alexander polynomials induced by periodic representations}
\vspace{-5pt}
\emph{Throughout this section, $P\co  \pi \to S_N$ is assumed to be a permutation representation induced by 
a finite representation $p\co  \pi' \to \Sigma$ of period $r$, as in \fullref{summary}.}

The representation $P$ induces an action of $\pi$ on the standard basis ${\cal B} = \{e_1, \ldots, e_N\}$ for $V={\mathbb Z}^N$. We obtain a representation $\rho\co  \pi \to GL(V)$. Let $\epsilon$ be  the abelianization homomorphism $\pi \to \<t\mid \>$  mapping $x \mapsto t$.  A twisted chain complex $C_*(X; \L \otimes V)$ is defined as in \fullref{review}. 

The free ${\mathbb Z}[\pi]$--complex $C_*(\tilde X)$ has a basis $\{\tilde z\}$ consisting of a  single lift of each cell $z$ in $X$.  Then $\{1 \otimes e_i \otimes \tilde z\}$ is a basis for the free 
$\L$--complex $C_*(X; \L\otimes V)$ (cf page 640 of \cite{kirkliv}).

We will use the following lemma from \cite{lift}. 

\begin{lemma}\label{module} Suppose that $A$ is a finitely generated ${\mathbb Z}[t^{\pm 1}]$--module admitting a square presentation matrix and has $0$th characteristic polynomial $\Delta(t)= c_0 \prod (t-\a_j)$. Let $s= t^r$, for some positive integer $r$. Then the $0$th characteristic polynomial of $A$, regarded  as a ${\mathbb Z}[s^{\pm 1}]$--module, is $\tilde \Delta(s)=
c_0^r \prod (s - \a_j^r)$. \end{lemma}

The map $P \co  \pi \to S_N$  restricts to a representation of the fundamental group $\pi'$ of the universal abelian cover $X_\infty$. Let $\hat X_\infty$ denote the induced $N$--fold cover. The $\L$--modules $H_1(\hat X_\infty)$ and 
$H_1(X; \L\otimes V)$ are isomorphic by two applications of 
Shapiro's Lemma (see for example \cite{hln}).

\begin{prop}\label{squarepres} $H_1(\hat X_\infty)$ is a finitely generated ${\mathbb Z}[s^{\pm 1}]$--module with a square presentation matrix, where $s= t^r$.
\end{prop}

\begin{proof} Construct $X_\infty$ in the standard way, splitting 
$X$ along the interior of Seifert surface  $S$ to obtain a relative cobordism $(V; S', S'')$ bounding two copies $S', S''$ of $S$.
Then $X_\infty$ is obtained by gluing 
countably many copies $(V_j; S'_j, S''_j)$ end-to-end, identifying 
$S''_j$ with $S'_{j+1}$, for each $j \in {\mathbb Z}$. 

For each $j$, let $W_j =
V_{jr}\cup \cdots \cup V_{jr +r -1}$ be the submanifold of $X_\infty$ bounding $S'_{jr}$ and $S''_{jr+r-1}$. 
Then $X_\infty$ is the union of the $W_j$'s, glued end-to-end. After lifting powers of the meridian of $k$, thereby constructing basepaths from $S_0'$ to each $S_{jr}' \subset W_j$, we can then regard each $\pi_1W_j$ as a subgroup of $\pi_1X_\infty \cong \pi'$. 

Conjugation by $x$ in the knot group induces an automorphism of $\pi'$, and the $r$th power maps $\pi _1W_j$ isomorphically to $\pi_1 W_{j+1}$. Since $p$ has period $r$, we have $p(x^{-r}ux^r)= p(u)$ for all $u \in \pi'$. Hence $P$ has the same image on each $\pi_1W_j$. 
By performing equivariant ambient $0$--surgery in $W_j$ to the lifted surfaces $\hat S_j'$ (that is, adding appropriate hollow 1--handles to the surface), we can assume that the image of $P (\pi_1 S_j')$ acts transitively, and hence each preimage $\hat S_j' \subset \hat X_\infty$ is connected. 

The covering space $\hat X_\infty$ is the union of countably many copies $\hat W_j$ of the lift $\hat W_0$ glued end-to-end. The cobordism $\hat W_0$, which bounds two copies $\hat S', \hat S''$ of  the surface $\hat S$, can be constructed from $\hat S' \times I$ by attaching $1$-- and $2$--handles in equal numbers. Consequently, $H_1 \hat W_0$ is a finitely generated  abelian group with a presentation of deficiency $d$ (number of generators minus number of relators) equal to the rank of $H_1 \hat S'$. 

The $r$th powers of covering transformations of $\hat X_\infty$ induce a ${\mathbb Z}[s^{\pm 1}]$--module structure on $H_1 \hat X_\infty$. The Mayer--Vietoris theorem implies that the generators of  $H_1 \hat W_0$ serve as generators for the module. Moreover, the relations of $H_1 \hat W_0$ together with $d$ relations arising from the boundary identifications become an equal number of relators. 
\end{proof}

\begin{cor}\label{alexone} If $\Delta_\rho(t) = 1$, then $H_1(\hat X_\infty)$ is trivial. \end{cor}

\begin{proof} Let $s = t^r$, and regard $H_1(\hat X_\infty)$ as a ${\mathbb Z}[s^{\pm 1}]$--module. Since the module has a square presentation matrix, its 
order ideal is principal, generated by $\tilde \Delta_\rho(s)$. \fullref{module} implies that $\tilde \Delta_\rho(s) =1$. Hence the order ideal coincides with the coefficient ring ${\mathbb Z}[s, s^{-1}]$. However, 
the order ideal is contained in the annihilator of the module (see \cite{crowell2} or
Theorem 3.1 of \cite{hillman}). Thus $H_1(\hat X_\infty)$ is trivial. 
\end{proof}

Since $p$ factors through $\pi_1 M_r$, so does $P|_{\pi'}$.  Let $\hat M_r $ denote the corresponding $N$--fold cover.

\begin{lemma}\label{quotient} $H_1\hat M_r$ is a quotient of $H_1 \hat X_\infty/(t^r-1)H_1 \hat X_\infty.$ \end{lemma}

\begin{proof}  Recall that $\pi_1 M_r \cong \pi'/[\pi', x^r]$.  Thus $\pi_1\hat M_r \cong \ker (P\vert_{\pi'})/[\pi',x^r]$, and by the Hurewicz theorem, $$H_1 \hat M_r \cong \ker(P\vert_{\pi'})/ \ker(P\vert_{\pi'})'\cdot [\pi', x^r].$$
On the other hand, $\pi_1\hat X_\infty$ modulo the relations $x^{-r} g x^r = g$ for all $g \in \pi_1\hat X_\infty$ is isomorphic to $ \ker (P\vert_{\pi'})/[ \ker (P\vert_{\pi'}),x^r]$.  Using the Hurewicz theorem again,
$$H_1\hat X_\infty/(t^r-1)H_1\hat X_\infty \cong \ker(P\vert_{\pi'})/ 
\ker(P\vert_{\pi'})'\cdot [\ker(P\vert_{\pi'}), x^r].$$
The conclusion follows immediately. 
\end{proof}

\begin{example} The group $\pi$ of the trefoil has presentation $\<x, a \mid  ax^2 a = xax\>$, where $x$ represents a meridian, and $a$ is in the commutator subgroup $\pi'$. The Reidemeister--Schreier method yields the presentation
$$\pi' = \<a_j\mid a_ja_{j+2}=a_{j+1}\>,$$
where $a_j=x^{-j}ax^j$.
Consider the homomorphism $p\co  \pi'\to \Sigma= \<\alpha \mid \alpha^3\> \cong {\mathbb Z}_3$ sending $a_{2j}\mapsto \a$ and $a_{2j+1}\mapsto \a^2$. We extend $p$ to $P \co  \pi \to  \hat \Sigma = \<\xi \mid \xi^2\> \ltimes \Sigma^2$, sending $x \mapsto \xi$.  The image $P (\pi')$  consists of the three elements $(1,1), (\a, \a^2), (\a^2, \a)$; the image of $\pi$ is isomorphic to the dihedral group $D_3$, which we regard as a subgroup of $GL_3({\mathbb Z})$. Hence we have a representation $\rho\co  \pi  \to GL_3({\mathbb Z})$. Let $\epsilon\co  \pi \to \<t \mid \>$ be the abelianization homomorphism mapping $x \mapsto t$. The product of $\rho$ and $\epsilon$ determines a tensor representation $\rho\otimes \epsilon\co  \pi \to GL_3({\mathbb Z}[t^{\pm 1}])$ defined by $(\rho\otimes \epsilon)(g) = \rho(g)\epsilon(g)$, for $g \in \pi$.  We order our basis so that: 
$$(\rho\otimes \epsilon)(x)=\begin{pmatrix} 0 & t & 0\\ t & 0 & 0\\ 0 & 0 & t\end{pmatrix},  \quad (\rho\otimes \epsilon)(a)=\begin{pmatrix} 0 & 1 & 0\\ 0 & 0 & 1\\ 1 & 0 & 0\end{pmatrix}$$
We can assume that the CW structure on $X$ contains a single $0$--cell $p$, $1$--cells $x, a$ and a single $2$--cell $r$.

The $\rho$--twisted cellular chain complex $C_*(X; \L\otimes V)$ has the form
$$0 \to C_2\cong \L^3\ {\buildrel \partial_2 \over \longrightarrow}\ C_1\cong \L^6 \ {\buildrel \partial_1\over \longrightarrow}\ C_0\cong \L^3 \to 0.$$ If we treat elements of $\L^3$ and $\L^6$ as row vectors, then
the map $\partial_2$ is described by a $3\times 6$  matrix obtained in the usual way from the $1\times 2$ matrix of Fox free derivatives:
$$( \frac {\partial r} { \partial x} \quad  \frac {\partial r}{ \partial a})= (a + ax -1 -xa \quad 1 + a x^2 -x)$$
replacing $x, a$ respectively with  their images under 
$\rho\otimes \epsilon$. The result is:
$$\partial_2 = \begin{pmatrix}  t-1 & 1 & -t &1 & t^2-t & 0  \\ 
 0 & -t-1 & t+1&-t & 1 & t^2 \\
1-t & t & -1&t^2 & 0 & 1-t  \end{pmatrix}$$ 
The map $\partial_1$ is determined by $(\rho\otimes \epsilon)(x) -I$ and $(\rho\otimes \epsilon)(a)- I$:  
$$\partial_1 = \begin{pmatrix}
-1&t & 0\\
t & -1 & 0\\
0 & 0 & t-1 \\
-1 & 1 & 0 \\
0 & -1 & 1 \\
1 & 0 & -1 \end{pmatrix}$$
Dropping the first three columns of the matrix for $\partial_2$ produces a $3\times 3$ matrix:
$$A= \begin{pmatrix} 1 & t^2 -t & 0 \\ 
-t & 1 & t^2 \\
t^2 & 0 & 1-t\end{pmatrix}$$
Similarly, eliminating the last three rows of $\partial_1$ gives:
$$B=\begin{pmatrix} -1 & t & 0 \\ 
t& -1 & 0 \\
0 & 0 & t-1\end{pmatrix}$$
Theorem 4.1 of \cite{kirkliv} implies that 
$\Delta_\rho(t)/ \Delta_0(t) = {\rm Det}\ A / {\rm Det}\ B$,
where $\Delta_0(t)$ is the $0$th characteristic polynomial of 
$H_0 \hat X_\infty$. Since $\hat X_\infty$ is connected, $\Delta_0(t) = t-1$. Hence $\Delta_\rho(t) = (t^2-t+1)(t^2-1).$

In this example, the cyclic resultant ${\rm Res}(\Delta_\rho(t), t^2-1)$ vanishes, indicating that 
$H_1\hat X_2$ is infinite. A direct calculation reveals that in fact $H_1 \hat X_2 \cong {\mathbb Z} \oplus {\mathbb Z}$.

 \end{example} 

\begin{remark}
In the above example we see that the Alexander polynomial of the trefoil knot divides the twisted Alexander polynomial. Generally, the Alexander polynomial divides any twisted Alexander polynomial arising from a finite permutation representation of the knot group.  A standard argument using the transfer homomorphism  and the fact that $H_1 X_\infty$ has no ${\mathbb Z}$--torsion shows that
$H_1 X_\infty$ embeds as a submodule in $H_1 \hat X_\infty$. Hence $\Delta(t)$, which is the 0th characteristic polynomial of $H_1 X_\infty$, divides $\Delta_\rho(t)$, the 0th characteristic polynomial
of $H_1\hat X_\infty$. \end{remark}

\section{Proof of the Theorem}

Alexander polynomials are a special case of twisted Alexander polynomials corresponding to the trivial representation. Hence it suffices to consider an arbitrary nontrivial knot $k$ with unit Alexander polynomial $\Delta(t)$.  

A complete list of those finite groups that can act freely on  a homology $3$--sphere is given in \cite{milnor}. The only nontrivial such group that is perfect (that is, has trivial abelianization) is the binary icosahedral group $A^*_5$, with order 120. 

Since $\Delta(t)$ annihilates $H_1 X_\infty$,  the condition that $\Delta(t)=1$ implies that $H_1 X_\infty$ is trivial or equivalently, that $\pi'$ is perfect. Hence each branched cover $M_r$ has perfect fundamental group and so is a homology sphere. 
Theorem 3.7 of \cite{lackenby} implies that for some integer $r > 2$, the group $\pi_1M_r$ is ``large" in the sense that it contains a finite-index subgroup with a free nonabelian quotient. 

Any large group has normal subgroups of arbitrarily large finite index.  Hence $\pi_1 M_r$ contains a normal subgroup $Q$ of index $N_0$ exceeding 120. Composing  the canonical projection $\pi' \to \pi_1 M_r$ with the quotient map $ \pi_1 M_r\to \pi_1 M_r/Q = \Sigma$, we obtain a surjective homomorphism $p\co \pi' \to \Sigma$ of least period $r_0$ dividing $r$. By \fullref{period1}, we have $r_0>1$. Let $P $ be the extension to $\pi$, as in \fullref{summary}. 
By that lemma, the order $N$ of $P (\pi')$ is no less than $N_0=| p(\pi')|$.  

As in section 2, realize $P (\pi)$ as a group of permutation matrices in $GL_N({\mathbb Z})$ acting transitively on the standard basis of ${\mathbb Z}^N$. Let $\hat M_{r_0}$ be the cover of $M_{r_0}$ induced by the representation $P \co  \pi \to S_N$. The group of covering transformations  acts freely on $\hat M_{r_0}$ and transitively on any point-preimage of the projection $\hat M_{r_0} \to M_{r_0}$. Its cardinality is equal $N$ and so cannot be the binary icosahedral group. Hence  $\hat M_{r_0}$ has nontrivial homology. 

\fullref{quotient} and \fullref{alexone} imply that $\Delta_\rho(t) \ne 1.$ \qed

\bibliographystyle{gtart}
\bibliography{link}

\end{document}